# Faster FISTA *

1st Jingwei Liang
*DAMTP, University of Cambridge*
Cambridge, United Kingdom
jl993@cam.ac.uk

2nd Carola-Bibiane Schönlieb
*DAMTP, University of Cambridge*
Cambridge, United Kingdom
cbs31@cam.ac.uk

***Abstract*—The "fast iterative shrinkage-thresholding algorithm", a.k.a. FISTA, is one of the most widely used algorithms in the literature. However, despite its optimal theoretical $O(1/k^2)$ convergence rate guarantee, oftentimes in practice its performance is not as desired owing to the (local) oscillatory behaviour. Over the years, various approaches are proposed to overcome this drawback of FISTA, in this paper, we propose a simple yet effective modification to the algorithm which has two advantages: 1) it enables us to prove the convergence of the generated sequence; 2) it shows superior practical performance compared to the original FISTA. Numerical experiments are presented to illustrate the superior performance of the proposed algorithm.**

## I. Introduction

### A. Problem statement

In various fields of science and engineering, including signal/image processing, inverse problems and machine learning, many problems to be solved can be cast as a structured composite non-smooth optimization problem of the sum of two functions, which usually reads

$$\min_{x\in\mathcal{H}} \ \Phi(x) \stackrel{\text{def}}{=} R(x) + F(x), \qquad (\mathcal{P})$$

where $\mathcal{H}$ is a real Hilbert space,

**(A.1)** $R : \mathcal{H} \to \mathbb{R} \cup \{+\infty\}$ is proper closed and convex;
**(A.2)** $F : \mathcal{H} \to \mathbb{R}$, and $\nabla F$ is $L$-Lipschitz continuous;
**(A.3)** $\mathrm{Argmin}(\Phi) \neq \emptyset$, the set of minimizers is non-empty.

Typical examples of $(\mathcal{P})$ can be found in Section IV.

### B. Forward–Backward splitting and FISTA schemes

One classical approach to solve problem $(\mathcal{P})$ is the Forward–Backward splitting (FBS) method [7], whose non-relaxed iteration takes the form

$$x_{k+1} \stackrel{\text{def}}{=} \mathrm{prox}_{\gamma_k R}\big(x_k - \gamma_k \nabla F(x_k)\big), \ \gamma_k \in ]0, 2/L[, \qquad (1)$$

where $\gamma_k$ is the step-size, and $\mathrm{prox}_{\gamma R}(\cdot) \stackrel{\text{def}}{=} \min_{x\in\mathcal{H}} \frac{1}{2}\|x - \cdot\|^2 + \gamma R(x)$ denotes the *proximity operator* of $R$.

The convergence of FB iterates is guaranteed as long as $\gamma_k \in [\epsilon_1, 2/L - \epsilon_2]$, $\epsilon_1, \epsilon_2 \in ]0, 1/L[$ is satisfied. For the convergence rate, it is well established that the objective function value, *i.e.* $\Phi(x_k) - \inf_{x\in\mathcal{H}} \Phi(x)$, converges at the speed of $O(1/k)$ which is quite slow. Over the years, various schemes have been proposed to accelerate the method. Among them, the FISTA scheme [3] by Beck and Teboulle (henceforth denoted as "**FISTA-BT**") is the most-known one, which achieves $O(1/k^2)$ convergence rate for the objective function value. However, while FBS is sequence convergent, the convergence property of the sequence generated by FISTA-BT has been a long-standing open problem. This question was settled in [4], followed by [2] in the continuous dynamical system case.

**Algorithm 1:** FISTA-BT algorithm

**Initial**: $t_0 = 1$, $\gamma = 1/L$ and $x_0 \in \mathcal{H}, x_{-1} = x_0$.
**repeat**

$$\begin{aligned} t_k &= \frac{1+\sqrt{1+4t_{k-1}^2}}{2}, \ a_k = \frac{t_{k-1}-1}{t_k}, \\ y_k &= x_k + a_k(x_k - x_{k-1}), \\ x_{k+1} &= \mathrm{prox}_{\gamma R}\big(y_k - \gamma\nabla F(y_k)\big). \end{aligned} \qquad (2)$$

$k = k + 1$;
**until** *convergence*;

In [4], Chambolle and Dossal proposed using the following strategy for the updating of $a_k$, let $d > 2$,

$$t_k = \frac{k+d}{d}, \ a_k = \frac{t_{k-1}-1}{t_k}. \qquad (3)$$

The above choice is denoted as "**FISTA-CD**". Under this setting, they managed to prove the convergence of the sequence while maintaining the $O(1/k^2)$ rate on the objective function values. Later in [1], the rate is proven to be actually $o(1/k^2)$.

### C. Slow practical performance

In practice, it has been reported in several work [10], [11], that despite the $O(1/k^2)$ convergence rate guarantee, oftentimes FISTA-BT has very slow practical performance which is mainly caused by the oscillation behaviour of the scheme. For the FISTA-CD scheme, when $d$ is close to 2, it has almost the same performance as FISTA-BT; see Section IV the numerical experiments.

In [6], it is reported that when $d$ is chosen in a certain range, such as $[50, 100]$, then practically FISTA-CD can be much faster than FISTA-BT; see also Section IV. A natural question would be raised: is it possible that the original FISTA-BT method can also achieve such boost of performance in

This work was partly supported by Leverhulme Trust project "Breaking the non-convexity barrier", the EPSRC grant "EP/M00483X/1", EPSRC centre "EP/N014588/1", the Cantab Capital Institute for the Mathematics of Information, and the Global Alliance project "Statistical and Mathematical Theory of Imaging".

practice? The main purpose of the presented paper is to answer this question.

## II. A MODIFIED FISTA SCHEME

In this section, we present the main contribution of this paper, a modified scheme of FISTA-BT.

### A. Two observations

In the original FISTA-BT scheme, the update of $t_k$ reads

$$t_k = \frac{1+\sqrt{1+4t_{k-1}^2}}{2}.$$

We have the following observations by replacing the $1, 1, 4$ in the numerator with parameters $p, q, r$.

*a) Parameter $r$:* Let $r > 0$, and consider $t_k = (1 + \sqrt{1+rt_{k-1}^2})/2$, then

$$a_k = \frac{t_{k-1}-1}{t_k} \begin{cases} r \in ]0,4[ : t_k \to \frac{4}{4-r},\ a_k \to \frac{r}{4}, \\ r \in [4,+\infty[ : t_k \to +\infty,\ a_k \to \frac{2}{\sqrt{r}}. \end{cases} \quad (4)$$

**Observation I:** $r$ controls the limiting value of $a_k$.

*b) Parameter $p, q$:* Let $r \in ]0,4]$ and $p, q > 0$. Consider $t_k = (p + \sqrt{q+rt_{k-1}^2})/2$, then

$$a_k = \frac{t_{k-1}-1}{t_k} \begin{cases} r \in ]0,4[ : t_k \to \frac{2p+\Delta}{4-r},\ a_k \to 1 - \frac{4-r}{2p+\Delta}, \\ r = 4 : t_k \to +\infty,\ a_k \to 1, \end{cases} \quad (5)$$

where $\Delta = \sqrt{rp^2 + (4-r)q}$.

Fix $r = 4$, Figure 1 shows the effects of different values of $p, q$, together with two different choices of $d$ in (3). Since $r = 4$, we have $a_k \to 1$, clearly, the *smaller* the value of $p, q$, the *slower* $a_k$ converges to 1. While for FISTA-CD, the *bigger* the value of $d$, the *slower* the $a_k$ converges to 1.

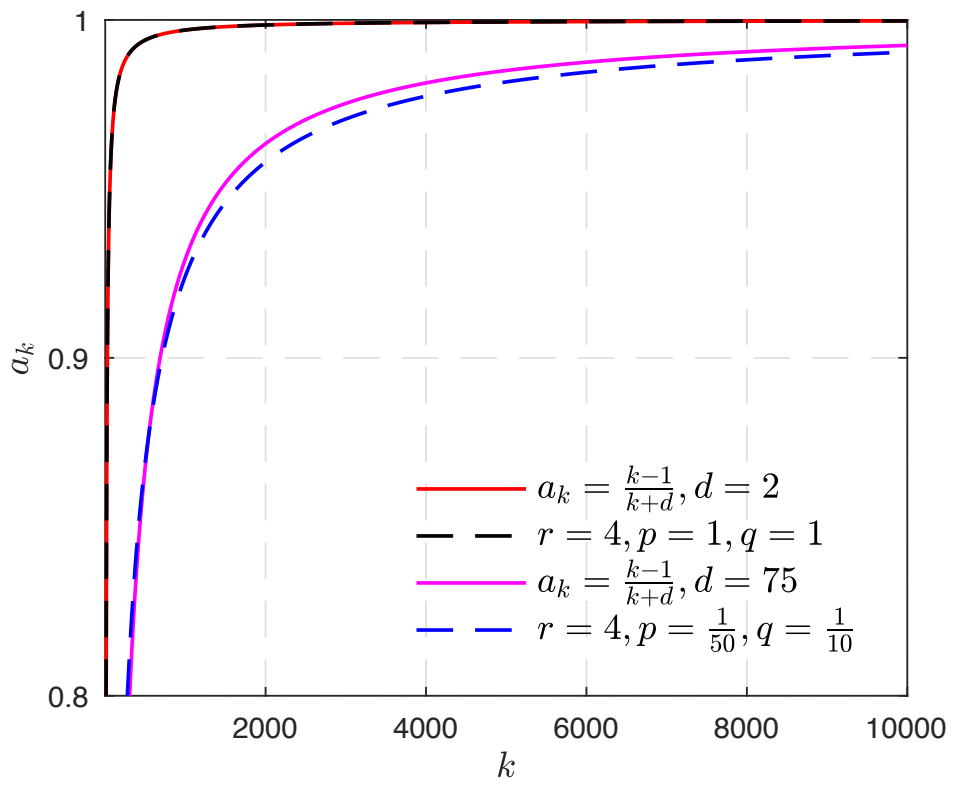

Fig. 1: Different effects of $p, q$ and $d$.

**Observation II:** $p, q$ control the speed of $a_k$ converging to 1.

**Remark II.1.** Fix $r = 4$, and if $q \le (2-p)^2$, then,

$$t_k = \frac{p+\sqrt{q+4t_{k-1}^2}}{2} \implies t_k^2 - t_k \le t_{k-1}^2,$$

which is the key to prove the convergence of the sequence of the modified FISTA scheme (Algorithm 2).

### B. A modified FISTA scheme

Based on the above observations, we propose the a modified FISTA scheme which is described in Algorithm 2. From now on, to distinguish Algorithm 2 from FISTA-BT and FISTA-CD, we shall call it as "FISTA-Mod".

**Algorithm 2:** A modified FISTA scheme

**Initial**: $p \in ]0,1], q > 0$ and $r \in ]0,4]$, $t_0 = 1$, $\gamma = 1/L$ and $x_0 \in \mathcal{H}, x_{-1} = x_0$.

**repeat**

$$\begin{aligned} t_k &= \frac{p+\sqrt{q+rt_{k-1}^2}}{2},\quad a_k = \frac{t_{k-1}-1}{t_k}, \\ y_k &= x_k + a_k(x_k - x_{k-1}), \\ x_{k+1} &= \mathrm{prox}_{\gamma R}\big(y_k - \gamma\nabla F(y_k)\big). \end{aligned} \quad (6)$$

**until** *convergence*;

**Remark II.2.** When $r \in ]0,4[$, then Algorithm 2 becomes a variant of the inertial Forward–Backward method [6].

### C. Convergence rate of the objective function

In this part we present the global convergence properties of FISTA-Mod scheme. We first show that FISTA-Mod preserves the $O(1/k^2)$ optimal convergence rate of FISTA-BT, and prove the convergence of the sequence $\{x_k\}_{k\in\mathbb{N}}$.

**Theorem II.3 (Convergence rate of the objective).** *For the FISTA-Mod scheme* (6)*, let $r = 4$ and choose $p \in ]0,1], q \in ]0,(2-p)^2]$. Then*

$$\Phi(x_k) - \Phi(x^\star) \le \frac{2L}{p^2(k+1)^2}\|x_0 - x^\star\|^2. \quad (7)$$

**Remark II.4.** Compared to the original convergence rate of FISTA-BT [3], which is $\Phi(x_k) - \Phi(x^\star) \le \frac{2L}{(k+1)^2}\|x_0 - x^\star\|^2$. Parameter $p$ appears in the obtained rate estimation, and $p = 1$ yields the *smallest* constant in the rate. Though $p < 1$ gives bigger constant in the rate estimation, as we shall see below it allows us to prove the $o(1/k^2)$ convergence rate.

**Sketch of proof.** There two key conditions to prove Theorem II.3, which are:

- from $q \le (2-p)^2$, one can show that $t_k^2 - t_k \le t_{k-1}^2$;
- For $p \in ]0,1]$, we have $t_k \ge \frac{k+1}{2}p$.

With the above result and follow the proof of [3, Theorem 4.1] we can prove Theorem II.3. □

**Theorem II.5 (From $O(1/k^2)$ to $o(1/k^2)$).** *For the FISTA-Mod scheme* (6)*, let $r = 4$ and choose $p \in ]0,1[, q > 0$ such that $p^2 \le q$. Then $\Phi(x_k) - \Phi(x^\star) = o(1/k^2)$.*

**Sketch of proof.** The key to establish $o(1/k^2)$ convergence rate is that with $0 < p < 1$, one can show that

$$\frac{p(1-p)(k+1)}{2} \le (1-p)t_k \le t_{k-1}^2 - (t_k^2 - t_k).$$

Then following [1], [4] we obtain the desired result. □

**Remark II.6.** (i) For the original FISTA-BT scheme, we have strictly $0 = t_{k-1}^2 - (t_k^2 - t_k)$, hence unable to obtain $o(1/k^2)$ convergence rate.

(ii) One byproduct of Theorem II.5 is one can show that sequence $\{x_k\}_{k\in\mathbb{N}}$ is bounded.

### D. Convergence of the sequence

**Theorem II.7 (Convergence of the sequence).** *For the FISTA-Mod scheme, let* $r = 4$ *and choose* $p \in ]0, 1[, q > 0$ *such that* $p^2 \le q$. *Then*

(i) *there exists an* $x^\star \in \mathrm{Argmin}(\Phi)$ *to which the sequence* $\{x_k\}_{k\in\mathbb{N}}$ *generated by FISTA-Mod converges weakly;*
(ii) *We have* $\|x_k - x_{k-1}\| = o(1/k)$.

**Sketch of proof.** There two key conditions to prove Theorem II.7, which are:

- sequence $\{x_k\}_{k\in\mathbb{N}}$ is bounded;
- The inertial parameter $\{a_k\}_{k\in\mathbb{N}}$ can be uniformly bounded from above by another sequence $\{a'_k\}_{k\in\mathbb{N}}$;

With the above result and follow the proof of [4, Theorem 4.1] we can prove Theorem II.7. □

## III. Lazy start and adaptive strategy

Since the modified FISTA scheme (6) has three degrees of freedom compared to the original FISTA-BT, we can design strategies to make FISTA-Mod adaptive to the properties of the problems so that faster practical performance can be achieved.

### A. Lazy-start FISTA-Mod

A very well-known behaviour of FISTA schemes is that, when the minimisation problem $(\mathcal{P})$ is strongly convex, both the trajectories of $\{\|x_k - x_{k-1}\|\}_{k\in\mathbb{N}}$ and $\{\Phi(x_k) - \Phi(x^\star)\}_{k\in\mathbb{N}}$ will oscillate once $a_k$ is too close to 1. Such oscillation slows down the speed of the algorithm [10], [11], and eventually makes the scheme slower than the original FBS scheme [6].

In [6], it is reported that for the FISTA-CD scheme, when the value of $d$ is chosen relative big (*e.g.* $d \in [50, 100]$), the FISTA-CD scheme achieves a much faster practical performance; see also the numerical experiments in Section IV. An intuitive explanation for such phenomena is that it is the interplay of the following two properties:

(i) for the considered problems in [6], they are locally strongly convex at the minimiser;
(ii) Relative bigger value of $d$ slows down the speed of $a_k$ converging to 1, as we have seen in Figure 1.

The interactions of the above two properties make the algorithm achieve faster practical performance.

Since all the parameters of the original FISTA-BT scheme are fixed, hence has fixed speed of $a_k$ converging to 1. While for the FISTA-Mod scheme, we can adjust the values of $p, q$ so that we can control the speed of $a_k$ approaching 1. In practice, we found the following choices of $p, q$ work quite well, which we dubbed as "lazy-start FISTA-Mod":

**Lazy-start FISTA-Mod:** $p \in [\frac{1}{100}, \frac{1}{10}], q \in ]0, 1]$.

### B. Adaptive to local strong convexity

In practice, many problems encountered are not globally strongly convex. However, oftentimes when certain conditions are satisfied (*e.g.* see for instance [6]), the problems locally have a so-called *quadratic growth* around the minimiser (see [6, Proposition 12]). As a result, locally adaptive strategies can be applied to achieve the optimal convergence rate.

Let us fix $\gamma = 1/L$, and suppose that problem $(\mathcal{P})$ is $\alpha$-strongly convex for some $\alpha > 0$, then the optimal choice of $a_k$ should be, according to [6, Section 4.4]

$$a_k \equiv a^\star = (1 - \sqrt{\gamma\alpha})^2/(1 - \gamma\alpha).$$

Define the following function of $\alpha$

$$f(\alpha) = 4(1 - \sqrt{\gamma\alpha})^2/(1 - \gamma\alpha).$$

Suppose that the local strong convexity of $\Phi$ is gradually changing until reaching $\alpha$, we propose the following adaptive FISTA scheme to take advantage of this local condition.

**Algorithm 3:** Adaptive-FISTA

**Initial**: $p = 1, q = 1$ and $r_0 = 4$, $t_0 = 1$, $\gamma = 1/L$ and $x_0 \in \mathcal{H}, x_{-1} = x_0$.
**repeat**
Estimate the local strong convexity $\alpha_k$;

$$r_k = f(\alpha_k),\ \ t_k = \frac{p + \sqrt{q + r_k t_{k-1}^2}}{2},\ \ a_k = \frac{t_{k-1} - 1}{t_k},$$
$$y_k = x_k + a_k(x_k - x_{k-1}),$$
$$x_{k+1} = \mathrm{prox}_{\gamma R}\big(y_k - \gamma\nabla F(y_k)\big).$$

**until** *convergence*;

For the rest of the paper, we shall call the above adaptive scheme "**Ada-FISTA**" for short.

## IV. Numerical experiments

In this section, we present numerical experiments of problems arising from linear inverse problem and image/video processing to demonstrate the advantages of the FISTA-Mod and Ada-FISTA over the original FISTA-BT.

### A. Linear inverse problem

We present first the numerical experiments of linear inverse problems. Consider the following forward observation of a vector $x_{\mathrm{ob}} \in \mathbb{R}^n$

$$f = \mathcal{K}x_{\mathrm{ob}} + w, \tag{8}$$

where $f \in \mathbb{R}^m$ is the observation, $\mathcal{K} : \mathbb{R}^n \to \mathbb{R}^m$ is some linear operator, and $w \in \mathbb{R}^m$ stands for noise. To recover or approximate $x_{\mathrm{ob}}$, one can consider the following optimization problem

$$\min_{x\in\mathbb{R}^n} \tfrac{1}{2}\|f - \mathcal{K}x\|^2 + \lambda R(x), \qquad (\mathcal{P}_\lambda)$$

where $\lambda > 0$ is the trade-off parameter, $R$ is the regulariser based on the prior knowledge on $x_{\mathrm{ob}}$.

We consider solving $(\mathcal{P}_\lambda)$ with $R$ being $\ell_1, \ell_{1,2}$-norms and $\ell_\infty$-norm. The observations are generated according to (8). Here $\mathcal{K}$ is generated from the standard Gaussian ensemble and the following parameters:

$\ell_1$**-norm** $(m, n) = (768, 2048)$, $x_{\mathrm{ob}}$ is 128-sparse;
$\ell_{1,2}$**-norm** $(m, n) = (512, 2048)$, $x_{\mathrm{ob}}$ has 16 non-zero blocks of size 8;

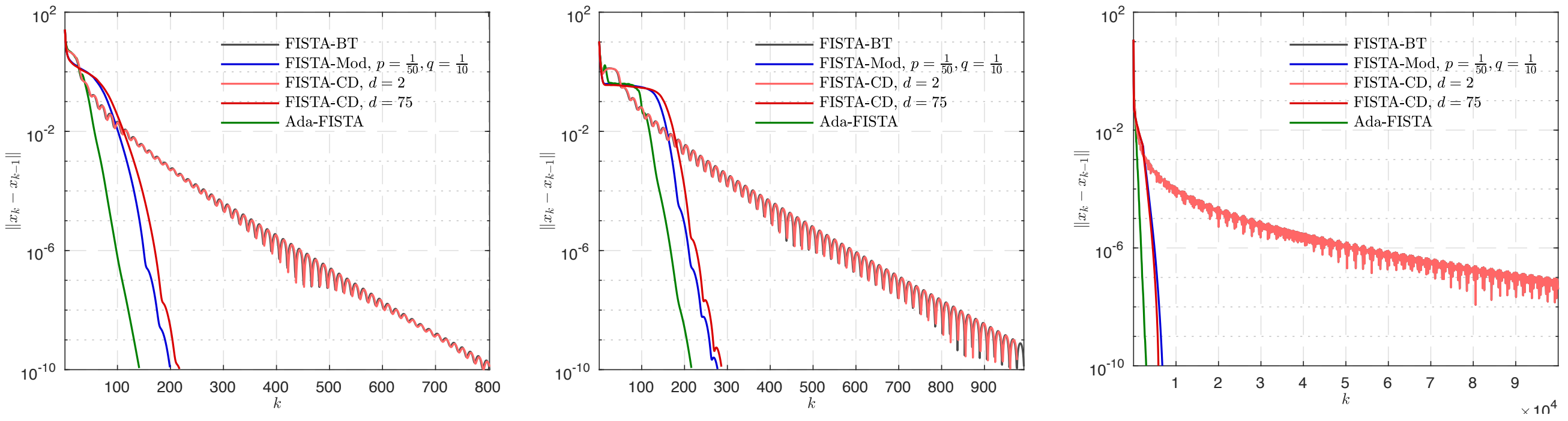

(a) $\ell_1$-norm (b) $\ell_{1,2}$-norm (c) $\ell_\infty$-norm

Fig. 2: Performance comparison of different FISTA schemes in terms of $\{\|x_k - x_{k-1}\|\}_{k\in\mathbb{N}}$ for linear inverse problem: (a) $\ell_1$-norm; (b) $\ell_{1,2}$-norm; (c) $\ell_\infty$-norm. The original FISTA-BT [3], sequence convergent FISTA-CD [4], and the proposed FISTA-Mod scheme. For FISTA-CD, two choices of $d$ are considered: $d = 2, 50$. For FISTA-Mod, $(p, q) = (1/50, 1/10)$ is considered. Black line, observation of FISTA-BT, the blue line is the observation of FISTA-Mod, the red lines are the observations of FISTA-CD, the green line is the observation of Ada-FISTA.

**$\ell_\infty$-norm** $(m, n) = (1020, 1024)$, $x_{\text{ob}}$ has 10 saturated entries.

The results of these examples are presented in Figure 2, from which we obtained the following observation and conclusions

- For all three examples, all the FISTA schemes exhibit a local linear convergence property, this is mainly due to the fact that considered $\ell_1, \ell_{1,2}, \ell_\infty$-norms belong to the so-called "partly smooth function", we refer to [6] for the dedicated study of this local linear convergence behaviour;
- The FISTA-CD with $d = 2$ has almost the same performance as the original FISTA-BT scheme, see the light red line and black line in all three figures;
- The FISTA-CD with $d = 75$ and the FISTA-Mod have very close performance, and both of them are much faster than FISTA-CD with $d = 2$ and the original FISTA-BT. More precisely, for $\ell_1, \ell_{1,2}$-norms, FISTA-CD with $d = 75$ and the FISTA-Mod are about 3 times faster, while for the $\ell_\infty$-norm, the difference is about 20 times which is quite significant;
- Ada-FISTA shows the fastest performance, especially for the $\ell_\infty$-norm, which is almost **50** times faster than the original FISTA-BT.

**Remark IV.1.** It should be noted that, a drawback of Ada-FISTA is that when the problem is of very large scale, estimating $\alpha_k$ each step can be very time consuming. A proper approach to deal with this deficiency is performing the evaluation in every $\kappa$ steps where $\kappa$ is properly chosen. For instance, for the experiments provided in Figure 2, $\alpha$ is estimated every 30 steps for $\ell_1, \ell_{1,2}$-norm and every 300 steps for $\ell_\infty$-norm.

### *B. Total variation based image deconvolution*

We also consider a 2D image processing problem, where $y$ is a degraded image generated according to (8), $\mathcal{K}$ is a circular convolution matrix with a Gaussian kernel. The anisotropic total variation (TV) [8] is applied for reconstruction, and the graph-cut algorithm [5] is applied for computing the proximity operator of TV.

The "cameraman" image is used for the experiments, the original, blurred and reconstructed images are shown in Figure 3(a)-(c). We compare only the performance of FISTA-BT and two settings of FISTA-Mod, the result is depicted in Figure 3(d). The result of this comparison is very similar to those of the linear inverse problem, the lazy-start FISTA-Mod shows superior performance than FISTA-BT.

### *C. Principal component pursuit*

To conclude this paper, we consider the principal component pursuit (PCP) problem [9], and apply it to decompose a video sequence into its background and foreground components.

Assume that a real matrix $y \in \mathbb{R}^{m\times n}$ can be written as

$$y = x_{\text{l,ob}} + x_{\text{s,ob}} + w,$$

where $x_{\text{l,ob}}$ is low–rank, $x_{\text{s,ob}}$ is sparse and $w$ is the noise. The PCP proposed in [9] attempts to provably recover $(x_{\text{l,ob}}, x_{\text{s,ob}})$ to a good approximation, by solving the following convex optimization problem

$$\min_{x_\text{l}, x_\text{s} \in \mathbb{R}^{m\times n}} \frac{1}{2}\|y - x_\text{l} - x_\text{s}\|_F^2 + \lambda_1\|x_\text{s}\|_1 + \lambda_2\|x_\text{l}\|_*, \tag{9}$$

where $\|\cdot\|_F$ is the Frobenius norm.

Observe that for fixed $x_\text{l}$, the minimizer of (9) is $x_\text{s}^\star = \text{prox}_{\lambda_1\|\cdot\|_1}(y - x_\text{l})$. Thus, (9) is equivalent to

$$\min_{x_\text{l}\in\mathbb{R}^{m\times n}} {}^1\big(\lambda_1\|\cdot\|_1\big)(y - x_\text{l}) + \lambda_2\|x_\text{l}\|_*, \tag{10}$$

where ${}^1\big(\lambda_1\|\cdot\|_1\big)(y - x_\text{l}) = \min_z \frac{1}{2}\|y - x_\text{l} - z\|_F^2 + \lambda_1\|z\|_1$ is the Moreau Envelope of $\lambda_1\|\cdot\|_1$ of index 1, and hence has 1-Lipschitz continuous gradient.

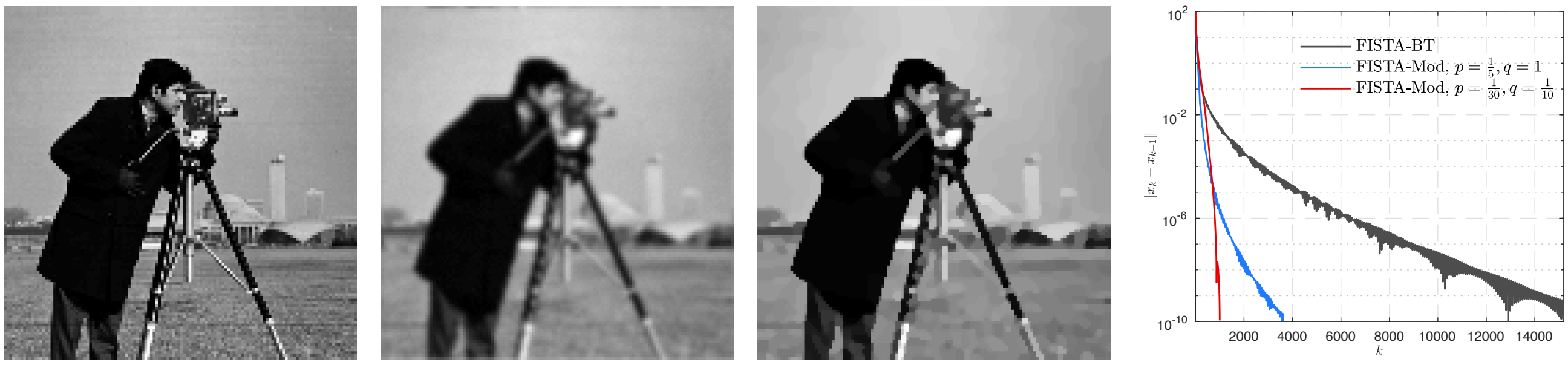

(a) Original image (b) Blurred image (c) Recovered image (d) Performance comparison

Fig. 3: Performance comparison of FISTA-BT and FISTA-Mod in terms of $\{\|x_k - x_{k-1}\|\}_{k\in\mathbb{N}}$ for TV based image deblurring: (a) original image; (b) blurred image; (c) deblurred image; (d) performance of FISTA-BT and FISTA-Mod.

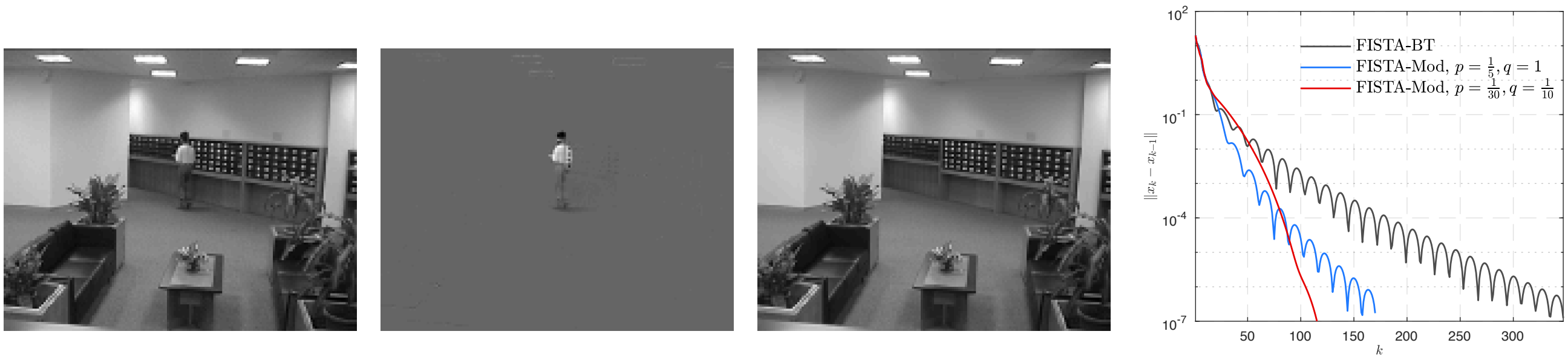

(a) Original frame (b) Sparse component (c) Low-rank component (d) Performance comparison

Fig. 4: Performance comparison of FISTA-BT and FISTA-Mod in terms of $\{\|x_k - x_{k-1}\|\}_{k\in\mathbb{N}}$ for principal component pursuit: (a) original frame; (b) sparse component; (c) low-rank component; (d) performance of FISTA-BT and FISTA-Mod.

We continue comparing only the performance of FISTA-BT and two settings of FISTA-Mod, the result is depicted in Figure 3(d). Again, the result of this comparison is very similar to previous examples, the lazy-start FISTA-Mod shows superior performance than FISTA-BT.